\documentclass[12pt]{article}
\usepackage{amssymb}
\usepackage{amsmath}

\newcommand{\lap}{\mbox{$\bigtriangleup$}}
\newcommand{\grad}{\mbox{$\bigtriangledown$}}
\newcommand{\ra}{{\mbox{$\rightarrow$}}}
\newcommand{\be}{\begin{equation}}
\newcommand{\ee}{\end{equation}}

\newtheorem{mthm}{Theorem}

\newtheorem{thm}{Theorem}[section]

\newtheorem{lem}{Lemma}[section]

\newtheorem{rem}{Remark}[section]
\begin{document}

\title{Maximum principles for the fractional p-Laplacian and symmetry of solutions}

\author{Wenxiong Chen \thanks{Partially supported by the Simons Foundation Collaboration Grant for Mathematicians 245486.} \quad  Congming Li \thanks{Corresponding author, School of Mathematics, Shanghai Jiao Tong University, congmingli@gmail.com, partially supported by NSFC 11571233 and NSF DMS-1405175.}}

\date{\today}
\maketitle
\begin{abstract}
In this paper, we consider nonlinear equations involving the fractional p-Laplacian
$$ (-\lap)_p^s u(x)) \equiv C_{n,s,p} PV \int_{\mathbb{R}^n} \frac{|u(x)-u(y)|^{p-2}[u(x)-u(y)]}{|x-z|^{n+ps}} dz= f(x,u).$$

We prove a {\em maximum principle for anti-symmetric functions} and obtain other key ingredients
for carrying on the method of moving planes, such as {\em a key boundary estimate lemma}. Then we establish radial symmetry and monotonicity for positive solutions to semilinear equations involving the fractional p-Laplacian in a unit ball and in the whole space.  We believe that the methods developed here
can be applied to a variety of problems involving nonlinear nonlocal operators.

\end{abstract}
\bigskip

\bigskip

\section{Introduction}

In this paper, we consider nonlinear equations involving the fractional p-Laplacian
\begin{equation}
(-\lap)_p^s u(x)) = f(x,u)
\label{Feq}
\end{equation}
with
\begin{eqnarray*}
 (-\lap)_p^s u(x)&=& C_{n,sp} \, \lim_{\epsilon \ra 0} \int_{\mathbb{R}^n\setminus B_{\epsilon}(x)} \frac{|u(x)-u(y)|^{p-2}[u(x)-u(y)]}{|x-y|^{n+sp}} dy \\
&=& C_{n,sp} PV \int_{\mathbb{R}^n} \frac{|u(x)-u(y)|^{p-2}[u(x)-u(y)]}{|x-y|^{n+sp}} dy,
\end{eqnarray*}
where PV stands for the Cauchy principal value.

In order the integral to make sense, we require that
$$u \in C^{1,1}_{loc} \cap L_{sp} $$
with
$$ L_{sp} = \{ u \in L^1_{loc} \mid \int_{\mathbb{R}^n} \frac{|u(x)|}{1+|x|^{n+sp}} d x < \infty \};$$

In the special case when $p=2$, $(-\lap)_p^s$ becomes the well-known fractional Laplacian $(-\lap)^s$. The nonlocal nature of these operators make them difficult to study. To circumvent this, Caffarelli and Silvestre \cite{CS} introduced the {\em extension method } which turns the nonlocal problem involving the fractional Laplacian into a local one in higher dimensions. This method has been applied successfully to study equations involving the fractional Laplacian, and a series of fruitful results have been obtained (see \cite{BCPS} \cite{CZ}    and the references therein). One can also use the {\em integral equations method}, such as {\em the method of moving planes in integral forms} and {\em regularity lifting} to investigate equations involving the fractional Laplacian by first showing that they are equivalent to the corresponding integral equations \cite{CLO} \cite{CLO1} \cite{CFY}.

However, when working at the extended problems or the corresponding integral equations, sometimes one needs to impose extra conditions on the solutions, which would not be necessary if one considers the pseudo differential equation directly (see the Introduction in \cite{CLL} for more details). Moreover, for equations involving the uniformly elliptic nonlocal operators
\begin{equation}
C_{n,\alpha} \, \lim_{\epsilon \ra 0} \int_{\mathbb{R}^n\setminus B_{\epsilon}(x)} \frac{a(x-z)((u(x)-u(z))}{|x-z|^{n+\alpha}} dz = f(x,u),
\label{a}
\end{equation}
where
$$ 0 < c_0 \leq a(y) \leq C_1, \; y \in R^n;$$
and for equations containing fully nonlinear nonlocal operators, such as
\begin{equation}
F_{\alpha}(u(x)) \equiv C_{n,\alpha} \, \lim_{\epsilon \ra 0} \int_{\mathbb{R}^n\setminus B_{\epsilon}(x)} \frac{G(u(x)-u(z))}{|x-z|^{n+\alpha}} dz =f(x,u)
\label{F}
\end{equation}
(see \cite{CS1} for the introductions of these operators), so far as we know, there has neither been any {\em extension method} nor {\em integral equation method} that work for these kinds of operators. This motivates us to come up with direct approaches on general nonlocal operators.

In our previous paper \cite{CLL}, a direct method of moving planes for the fractional Laplacian has been introduced and has been applied to obtain symmetry, monotonicity, and non-existence of solutions for various semi-linear equations involving the fractional Laplacian. Moreover, this direct approach can be applied to study the qualitative properties of solutions to uniformly elliptic problem (\ref{a}) \cite{TF}. In \cite{CLLg}, we refined this direct approach, so that it can be applied to fully nonlinear nonlocal problem (\ref{F}) in the case the operator is non-degenerate in the sense that
$$G'(t) \geq c_o > 0.$$
This includes the fractional Laplacian when $G(t)=t$ and $G'(t)=1$.

The fractional p-Laplacian we considered in this paper is actually a special case of the nonlinear nonlocal operator $F_{\alpha}(\cdot)$ when
$$\alpha = s p \; \mbox{ and } G(t) = |t|^{p-2}t .$$
However, it is degenerate for $p>2$ and singular for $p<2$. Here $G'(t) = (p-1)|t|^{p-2}$, and
$$ \mbox{ as } t \ra 0, \; G'(t) \ra \left\{\begin{array}{ll} 0 & \mbox{ if } p > 2 \\
\infty &  \mbox{ if } 1<p<2.
\end{array}
\right.
$$

The methods introduced in either \cite{CLL} or \cite{CLLg} relies heavily on the non-degeneracy of $G(\cdot)$, hence they cannot be applied directly to the degenerate or singular fractional p-Laplacian.
To overcome these difficulties, we introduce some new ideas, among which a significant one is the so-called {\em the key boundary estimate lemma}, and it plays the role of Hopf lemma in the second step of the method of moving planes. We believe that this kind of lemma and the idea behind the proof will become powerful tools in studying the qualitative properties of solutions for nonlinear nonlocal problems.

We first prove

\begin{mthm} (A Simple Maximum Principle)

Let $\Omega$ be a bounded domain in $\mathbb{R}^n$. Assume that $u \in C^{1,1}_{loc} (\Omega) \cap L_{sp}$, be lower semi-continuous on $\bar{\Omega}$, and satisfies
\begin{equation}
\left\{\begin{array}{ll}
(-\lap)^s_p u(x) \geq 0, & x \in \Omega,\\
u(x) \geq 0, & x \in \Omega^c.
\end{array}
\right.
\end{equation}

Then
\begin{equation}
u(x) \geq 0, \;\;\;\; x \in \Omega.
\end{equation}

If $u(x)=0$ at some point $x\in \Omega$, then
$$ u(x) = 0 \; \mbox{ almost everywhere in } \mathbb{R}^n.$$

The same conclusion holds for unbounded domains $\Omega$ if we further assume that
$$\underset{|x| \ra \infty}{\underline{\lim}} u(x) \geq 0.$$
\label{mthm1}
\end{mthm}

Then we establish {\em maximum principles} and {\em a key boundary estimate lemma} for anti-symmetric functions which play important roles in carrying out the method of moving planes. To explain the assumptions in these principles, we recall some basics in this method.

Take the whole space $\mathbb{R}^n$ as an example. Let
$$T_{\lambda} =\{x \in \mathbb{R}^{n}|\; x_1=\lambda, \mbox{ for some } \lambda\in \mathbb{R}\}$$
be the moving planes,
$$\Sigma_{\lambda} =\{x \in \mathbb{R}^{n} | \, x_1<\lambda\}$$
be the region to the left of the plane, and
$$ x^{\lambda} =(2\lambda-x_1, x_2, ..., x_n)$$
be the reflection of $x$ about the plane $T_{\lambda}$.

Assume that $u$ is a solution of pseudo differential equation (\ref{Feq}). To compare the values of $u(x)$ with $u(x^{\lambda})$, we denote
$$w_{\lambda} (x) = u(x^{\lambda}) - u(x) .$$
Obviously, $w_{\lambda}(x^{\lambda})=-w_{\lambda}(x)$, hence it is said to be anti-symmetric.

The first step is to show that for $\lambda$ sufficiently negative, we have
\begin{equation}
w_{\lambda}(x) \geq 0 , \;\; x \in \Sigma_{\lambda} .
\label{w}
\end{equation}
This provides a starting point to move the plane. Then in the second step, we move the plane to the right as long as inequality (\ref{w}) holds to its limiting position to show that $u$ is symmetric
about the limiting plane. To prove (\ref{w}), we will use

\begin{mthm} \label{mthm2} ( A Maximum Principle for Anti-symmetric Functions).

Let $\Omega$ be a bounded domain in $\Sigma_{\lambda}$.
Assume that $w\in L_{sp}\cap C^{1,1}(\Omega)$ and is lower semi-continuous on $\bar{\Omega}$. If
\be
\left\{\begin{array}{ll}
(-\lap)^s_p u_{\lambda}(x) - (-\lap)^s_p u(x) \geq0  &\mbox{ in } \Omega,\\
w_{\lambda}(x) \geq0&\mbox{ in } \Sigma_{\lambda} \backslash\Omega,
\end{array}
\right.
\label{s32}
\ee
then
\be
w_{\lambda}(x) \geq0 \mbox{ in } \Omega.\label{s41}
\ee

If $w_{\lambda} = 0$ at some point in $\Omega$, then
$$ w_\lambda(x) = 0 \; \mbox{almost everywhere in }  \mathbb{R}^n. $$

These conclusions hold for unbounded region $\Omega$ if we further assume that
$$\underset{|x| \ra \infty}{\underline{\lim}} w_\lambda(x) \geq 0.$$
\end{mthm}

In the second step, let
$$\lambda_o = \sup \{ \lambda \mid w_{\mu}(x) \geq 0, \; x \in \Sigma_{\mu}, \, \mu \leq \lambda\}$$
be the upper limit of such $\lambda$ that (\ref{w}) holds.
To show that $u$ is symmetric about the limiting plane $T_{\lambda_o}$, or
\begin{equation}
 w_{\lambda_o}(x) \equiv 0 , \;\; x \in \Sigma_{\lambda_o};
 \label{w0}
 \end{equation}
one usually uses a contradiction argument: if (\ref{w0}) does not hold, then one can move the
plane a little bit forward, and still have (\ref{w}) for some $\lambda > \lambda_o$. The region
between $T_{\lambda_o}$ and $T_{\lambda}$ is a narrow region, hence a {\em Narrow Region Principle} is often employed for this purpose. However, for the degenerate or singular fractional p-Laplacian, this kind of {\em Narrow Region Principle} can no longer be obtained, and one needs to carry out a finer analysis.

The main approach goes roughly as the following.

Suppose (\ref{w0}) is violated, then by a strong maximum principle, we have
$$w_{\lambda_o}(x) > 0 , \;\; x \in \Sigma_{\lambda_o}.$$

On the other hand, by the definition of $\lambda_o$, there exists a sequence $\lambda_k \searrow \lambda_o$, and $x^k \in \Sigma_{\lambda_k}$, such that
\begin{equation}
w_{\lambda_k}(x^k) = \min_{\Sigma_{\lambda_k}} w_{\lambda_k}  < 0, \;\; \mbox{ and } \; \grad w_{\lambda_k}(x^k) = 0 .
\label{wxk}
\end{equation}

By Theorem \ref{mthm2} and under a mild assumption on $f(x,u)$, we can derive that there is a subsequence of $\{x^k\}$ that converges to some point $x^o$, and from (\ref{wxk}), we have
$$w_{\lambda_o}(x^o) \leq 0, \; \mbox{ hence } x^o \in \partial \Sigma_{\lambda_o}; \;\; \mbox{ and } \; \grad w_{\lambda_o}(x^o) = 0 .$$
This will contradicts the following
\begin{mthm}(A key boundary estimate) \label{Hopf}

Assume that $w_{\lambda_o} > 0,$ for $x \in \Sigma_{\lambda_o}$. Suppose $\lambda_k \ra \lambda_o$, and $x^k \in \Sigma_{\lambda_k}$, such that
$$
w_{\lambda_k}(x^k) = \min_{\Sigma_{\lambda_k}} w_{\lambda_k}  \leq 0 \; \mbox{ and } x^k \ra x^o \in \partial \Sigma_{\lambda_o}.
$$

Let
$$\delta_k = dist(x^k, \partial \Sigma_{\lambda_k}) \equiv |\lambda_k -x^k_1|.$$

Then
$$
\overline{\lim_{\delta_k \ra 0}} \frac{1}{\delta_k} \left\{ (-\lap)^s_p u_{\lambda_k} (x^k) - (-\lap)^s_p u(x^k) \right\} < 0.
$$
\end{mthm}

We will use two examples to illustrate how these key ingredients in the above can be used in the {\em method of moving planes} to establish symmetry and monotonicity of positive solutions.

We first consider
\be
\left\{ \begin{array}{ll}
(-\lap)^s_p u(x)  = u^q(x) , & x \in B_1(0), \\
u(x) = 0 , & x \not{\in} B_1(0).
\end{array}
\right.
\label{mfu1}
\ee
We prove
\begin{mthm}
Assume that $u \in C_{loc}^{1,1}(B_1(0))$ is a positive solution of (\ref{mfu1}) with $q \geq p-1$. Then $u$ must be radially symmetric and monotone decreasing about the origin.
\label{mthm4}
\end{mthm}

Since the {\em narrow region principle} is no longer valid for the fractional p-Laplacian,
we employ some new ideas to get start the moving of the plane.

Then we study
\be
(-\lap)^s_p u(x) = g(u(x)), \;\; x \in \mathbb{R}^n .
\label{meqws}
\ee

\begin{mthm}
Assume that $u \in C^{1,1}_{loc} \cap L_{sp}$ is a positive solution of (\ref{meqws})with
$\lim_{|x| \ra \infty} u(x) = 0.$
Suppose $g'(s) \leq 0$ for $s$ sufficiently small.

Then $u$ must be radially symmetric and monotone decreasing about some point in $\mathbb{R}^n$.
\label{mthm5}
\end{mthm}

\begin{rem}
 In \cite{LN}, Li and Ni considered a similar semilinear equation for regular Laplacian
 $$- \lap u(x) = g(u(x)), \;\; x \in \mathbb{R}^n .$$
Under the same assumptions on $u$ and $g(\cdot)$, they obtained the radial symmetry for
positive solutions. Our theorem here generalize their result to the fractional p-Laplacian.
One can see that in the special case when $p=2$, and as $s \ra 1$, $(-\lap)^s u(x) \ra -\lap u(x)$.

\end{rem}

In Section 2, we establish various maximum principles and prove Theorem \ref{mthm1}, \ref{mthm2}, and \ref{Hopf}. Using the key ingredients obtained in Section 2, we derive symmetry and montonicity of positive solutions and prove Theorem \ref{mthm4} and \ref{mthm5} in Section 3 and 4 respectively.

For more articles concerning the method of moving planes for nonlocal equations and for integral equations, please see \cite{CLv} \cite{FL} \cite{Ha} \cite{HLZ} \cite{HWY} \cite{Lei} \cite{LL} \cite{LLM} \cite{LZ} \cite{LZ1} \cite{LZ2} \cite{MC} \cite{MZ} and the references therein.

\section{Various Maximum Principles and a Boundary Estimate}
Let
\begin{eqnarray}
 (-\lap)_p^s u(x))&=& C_{n,sp} \, \lim_{\epsilon \ra 0} \int_{\mathbb{R}^n\setminus B_{\epsilon}(x)} \frac{|u(x)-u(y)|^{p-2}[u(x)-u(y)]}{|x-y|^{n+sp}} dy \nonumber \\
&=& C_{n,sp} PV \int_{\mathbb{R}^n} \frac{|u(x)-u(y)|^{p-2}[u(x)-u(y)]}{|x-y|^{n+sp}} dy,
\label{2.1}
\end{eqnarray}

\begin{thm} (A Simple Maximum Principle)

Let $\Omega$ be a bounded domain in $\mathbb{R}^n$. Assume that $u \in C^{1,1}_{loc} (\Omega) \cap L_{sp}$ be lower semi-continuous on $\bar{\Omega}$, and satisfies
\begin{equation}
\left\{\begin{array}{ll}
(-\lap)^s_p u(x) \geq 0, & x \in \Omega,\\
u(x) \geq 0, & x \in \Omega^c.
\end{array}
\right.
\label{2.2}
\end{equation}

Then
\begin{equation}
u(x) \geq 0, \;\;\;\; x \in \Omega.
\label{2.3}
\end{equation}

If $u(x)=0$ at some point $x\in \Omega$, then
$$ u(x) = 0 \; \mbox{ almost everywhere in } \mathbb{R}^n.$$

The same conclusion holds for unbounded domains $\Omega$ if we further assume that
$$\underset{|x| \ra \infty}{\underline{\lim}} u(x) \geq 0.$$
\label{thm2.1}
\end{thm}

{\bf Proof.} Suppose (\ref{2.3}) is violated, then since $u$ is lower semi-continuous on $\bar{\Omega}$, there exists $x^o$ in $\Omega$ such that
$$u(x^o) = \min_{\Omega} u < 0.$$
It follows from the second inequality in (\ref{2.2}) that
\begin{eqnarray*}
(-\lap)^s_p u(x^o) &=& \int_{\mathbb{R}^n} \frac{|u(x^o)-u(z)|^{p-2}[u(x^o)-u(z)]}{|x^o-z|^{n+sp}} dz \\
&\leq& \int_{\Omega^c} \frac{|u(x^o)-u(z)|^{p-2}[u(x^o)-u(z)]}{|x^o-z|^{n+sp}} dz\\
&<& 0.
\end{eqnarray*}
This contradicts the first inequality in (\ref{2.2}) and hence (\ref{2.3}) must be true. It follows that if $u(x^o)=0$ at some point $x^o \in \Omega$, then
\begin{eqnarray*}
(-\lap)^s_p u(x^o) &=& \int_{\mathbb{R}^n} \frac{|u(z)|^{p-2}[-u(z)]}{|x^o-z|^{n+sp}} dz \\
&\leq& 0,
\end{eqnarray*}
while on the other hand, from the first inequality in (\ref{2.2}), we should have
$$\int_{\mathbb{R}^n} \frac{|u(z)|^{p-2}[-u(z)]}{|x^o-z|^{n+sp}} dz \geq 0,$$
and hence the integral must be $0$. Taking into account that $u$ is nonnegative, we derive
$$ u(x) = 0 \; \mbox{ almost everywhere in } \mathbb{R}^n.$$
This proves the theorem.

In the following, we will continue to use the notation introduced in the previous section. Let

$$T_{\lambda} =\{x \in \mathbb{R}^{n}|\; x_1=\lambda, \mbox{ for some } \lambda\in \mathbb{R}\}$$
be the moving planes,
$$\Sigma_{\lambda} =\{x \in \mathbb{R}^{n} | \, x_1<\lambda\}$$
be the region to the left of the plane, and
$$ x^{\lambda} =(2\lambda-x_1, x_2, ..., x_n)$$
be the reflection of $x$ about the plane $T_{\lambda}$.
$$w_{\lambda} (x) = u(x^{\lambda}) - u(x) .$$
For simplicity of notation, we denote $w_{\lambda}$ by $w$ and $\Sigma_{\lambda}$ by $\Sigma$.

\begin{thm}(Maximum Principle for Anti-symmetric Functions)
Let $\Omega$ be a bounded domain in $\Sigma$.
Assume that $w\in L_{sp}\cap C_{loc}^{1,1}(\Omega)$ and is lower semi-continuous on $\bar{\Omega}$. If
$$
\left\{\begin{array}{ll}
(-\lap)^s_p u_{\lambda}(x) - (-\lap)^s_p u(x)  \geq 0  &\mbox{ in } \Omega,\\
w(x) \geq0&  \mbox{ in }   \Sigma \backslash\Omega,\\
w(x^{\lambda})=-w(x)  &\mbox{ in } \Sigma,
\end{array}
\right.
$$
then
$$
w(x) \geq0 \mbox{ in } \Omega.
$$

Furthermore, if $w = 0$ at some point in $\Omega$, then
$$ w(x) = 0 \; \mbox{ almost everywhere in }  \mathbb{R}^n. $$

These conclusions hold for unbounded region $\Omega$ if we further assume that
$$\underset{|x| \ra \infty}{\underline{\lim}} w(x) \geq 0 .$$

\label{thm2.2}
\end{thm}

{\bf Proof.} Suppose otherwise, then there exists a point $x$ in $\Omega$, such that
\begin{equation}
w(x) = \min_{\Omega} w < 0 .
\label{w<0}
\end{equation}

For simplicity of writing, we denote $G(t)=|t|^{p-2}t$. Obviously, $G(t)$ is a strictly increasing function, and $G'(t) =(p-1)|t|^{p-2} \geq 0$.

\begin{eqnarray}
& & (-\lap)^s_p u_{\lambda}(x) - (-\lap)^s_p u(x) \nonumber \\
&=& C_{n, sp} PV \int_{\mathbb{R}^n} \frac{G(u_{\lambda}(x) -u_{\lambda}(y)) - G(u(x)-u(y))}{|x-y|^{n+sp}} d y \nonumber \\
&=& C_{n, sp} PV \int_{\Sigma} \frac{G(u_{\lambda}(x) -u_{\lambda}(y)) - G(u(x)-u(y))}{|x-y|^{n+sp}} d y\nonumber \\
 &+& C_{n, sp} PV \int_{\Sigma} \frac{G(u_{\lambda}(x) -u(y)) - G(u(x)-u_{\lambda}(y))}{|x-y^{\lambda}|^{n+sp}} d y \nonumber\\
&=& C_{n, sp} PV \int_{\Sigma} \left[ \frac{1}{|x-y|^{n+sp}} -  \frac{1}{|x-y^\lambda|^{n+sp}}\right] \left[G(u_{\lambda}(x) -u_{\lambda}(y)) - G(u(x)-u(y))\right]  d y\nonumber \\
&+& C_{n, sp} PV \int_{\Sigma} \frac{G(u_{\lambda}(x) -u_{\lambda}(y)) - G(u(x)-u(y))+ G(u_{\lambda}(x) -u(y)) - G(u(x)-u_{\lambda}(y))}{|x-y^{\lambda}|^{n+sp}} d y \nonumber\\
&=& C_{n, sp}  \left\{ I_1 + I_2 \right\}.
\label{2.5}
\end{eqnarray}

To estimate $I_1$, we first notice that
$$ \frac{1}{|x-y|} > \frac{1}{|x-y^{\lambda}|} , \;\; \forall \, x, y \in \Sigma.  $$
While for second part in the integral, we have
$$G(u_{\lambda}(x) -u_{\lambda}(y)) - G(u(x)-u(y)) \leq 0 \; \mbox{ but } \not{\equiv} 0;$$
due to the monotonicity of $G$ and the fact that
$$[u_{\lambda}(x) -u_{\lambda}(y)]-[u(x)-u(y)] = w(x)-w(y) \leq 0 \; \mbox{ but } \not{\equiv} 0.$$
It follows that
\begin{equation}
I_1 < 0 .
\label{I1}
\end{equation}

To estimate $I_2$, we regroup the terms and apply the {\em mean value theorem} to derive
\begin{eqnarray}
I_2 &=& \int_{\Sigma} \frac{[G(u_{\lambda}(x) -u_{\lambda}(y)) - G(u(x)-u_{\lambda}(y))]+[G(u_{\lambda}(x) -u(y))- G(u(x)-u(y))] }{|x-y^{\lambda}|^{n+sp}} d y \nonumber\\
&=& w (x) \int_{\Sigma} \frac{G'(\xi(y)) + G'(\eta(y))}{|x-y^{\lambda}|^{n+sp}} d y \leq 0.
\label{I2}
\end{eqnarray}

Combining (\ref{2.5}), (\ref{I1}), and (\ref{I2}), we deduce
$$ (-\lap)^s_p u_{\lambda}(x) - (-\lap)^s_p u(x) < 0 .$$
This contradicts our assumption, hence we must have
$$w(x) \geq 0.$$

It follows that if $w(x) =0$ at some point $x \in \Omega$, then $x$ is a minimum of $w$ in $\Omega$, hence (\ref{2.5}) holds with $I_2=0$. Now, our assumption implies $I_1 \geq 0$, consequently
$$G(u_{\lambda}(x) -u_{\lambda}(y)) - G(u(x)-u(y)) \geq 0,$$
and due to the monotonicity of $G$, we derive, for almost all $y \in \Sigma$,
$$[u_{\lambda}(x) -u_{\lambda}(y)]-[u(x)-u(y)] = w(x)-w(y) = -w(y) \geq 0.$$
Therefore, we must have
$$w (y) = 0 , \mbox{ almost everywhere in } \Sigma ,$$
and from the antisymmetry of $w_\lambda$, we arrive at
$$ w(y) = 0 \; \mbox{ almost everywhere in }  \mathbb{R}^n. $$
This completes the proof.

\begin{thm}(A Boundary Estimate) \label{thmHopf}

Assume that $w_{\lambda_o} > 0,$ for $x \in \Sigma_{\lambda_o}$. Suppose $\lambda_k \searrow \lambda_o$, and $x^k \in \Sigma_{\lambda_k}$, such that
\begin{equation}
w_{\lambda_k}(x^k) = \min_{\Sigma_{\lambda_k}} w_{\lambda_k}  \leq 0 \; \mbox{ and } x^k \ra x^o \in \partial \Sigma_{\lambda_o}.
\label{wk}
\end{equation}
Let
$$\delta_k = dist(x^k, \partial \Sigma_{\lambda_k}) \equiv |\lambda_k -x^k_1|.$$

Then
\begin{equation}
\overline{\lim_{\delta_k \ra 0}} \frac{1}{\delta_k} \left\{ (-\lap)^s_p u_{\lambda_k} (x^k) - (-\lap)^s_p u(x^k) \right\} < 0.
\label{hop}
\end{equation}
\end{thm}

\begin{rem} The limit in (\ref{hop}) may be $-\infty$.
\end{rem}

{\bf Proof.} Similar to the proof in Theorem \ref{thm2.2}, we derive

\begin{eqnarray}
& & \frac{1}{\delta_k} \left\{ (-\lap)^s_p u_{\lambda_k} (x^k) - (-\lap)^s_p u(x^k) \right\} \nonumber\\
&=& \frac{C_{n, sp}}{\delta_k} PV \int_{\Sigma_{\lambda_k}} \left[ \frac{1}{|x^k-y|^{n+sp}} -  \frac{1}{|x^k-y^{\lambda_k}|^{n+sp}}\right] \left[G(u_{\lambda_k}(x^k) -u_{\lambda_k}(y)) - G(u(x^k)-u(y))\right]  d y\nonumber \\
&+& \frac{w_{\lambda_k}(x^k)}{\delta_k} C_{n, sp} PV \int_{\Sigma_{\lambda_k}} \frac{G'(\xi(y)) + G'(\eta(y))}{|x^k-y^{\lambda_k}|^{n+sp}} d y \nonumber\\
&=& C_{n, sp} PV \left\{ I_{1k} + I_{2k} \right\}.
\label{2.6}
\end{eqnarray}

Obviously,
\begin{equation}
I_{2k} \leq 0 .
\label{2.7}
\end{equation}

To estimate $I_{1k}$, we apply the {\em mean value theorem} on $f(t)= t^{-(n+sp)/2}$ with $t=|\cdots|^2$ to obtain
\begin{eqnarray}
 & & \frac{1}{\delta_k} \left[\frac{1}{|x^k-y|^{n+sp}} -  \frac{1}{|x^k-y^{\lambda_k}|^{n+sp}}\right] \nonumber\\
 &=& \frac{2(n+sp)(\lambda_k-y_1)}{|x^k- \eta_k(y)|^{n+sp+2}} \ra \frac{2(n+sp)(\lambda_o-y_1)}{|x^o- \eta_o(y)|^{n+sp+2}}, \mbox{ as } k \ra \infty.
 \label{2.8}
\end{eqnarray}
Here $\eta_k(y)$ is some point on the line segment between $y$ and $y^{\lambda_k}$, and hence
$\eta_o(y)$ is some point on the line segment between $y$ and $y^{\lambda_o}$.

It is easy to see that the last term in (\ref{2.8}) is strict positive in $\Sigma_{\lambda_o}$ (it may be +$\infty$ at some point $y$).

Meanwhile, as $k \ra \infty$,
\begin{equation}
G(u_{\lambda_k}(x^k) -u_{\lambda_k}(y)) - G(u(x^k)-u(y)) \ra
G(u_{\lambda_o}(x^o) -u_{\lambda_o}(y)) - G(u(x^o)-u(y)) < 0
\label{2.9}
\end{equation}
for all $y \in \Sigma_{\lambda_o}$, because
$$\left[u_{\lambda_o}(x^o) -u_{\lambda_o}(y)\right]-\left[u(x^o)-u(y)\right] = w_{\lambda_o}(x^o)-w_{\lambda_o}(y) < 0,$$
and due to the monotonicity of $G(\cdot)$.

Combining (\ref{2.6}), (\ref{2.7}), (\ref{2.8}), and (\ref{2.9}), we arrive at (\ref{hop}).

This completes the proof of the theorem.

\section{Symmetry and Monotonicity in a Ball Region}

In this section, we consider
\be
\left\{ \begin{array}{ll}
(-\lap)^s_p u(x)  = u^q(x) , & x \in B_1(0), \\
u(x) = 0 , & x \not{\in} B_1(0).
\end{array}
\right.
\label{mfu1a}
\ee

\begin{thm}
Assume that $u \in C_{loc}^{1,1}(B_1(0))$ is a positive solution of (\ref{mfu1a}) with $q \geq p-1$. Then $u$ must be radially symmetric and monotone decreasing about the origin.
\label{thm4.1}
\end{thm}

{\bf Proof.} Let $\Sigma_{\lambda}$, $u_{\lambda}(x)$, and $w_{\lambda}(x)$ as defined in the previous section. Denote $$\Omega_{\lambda} = \Sigma_{\lambda} \cap B_1(0).$$

Then in $\Omega_\lambda$, we have
$$(-\lap)^s_p u_{\lambda}(x) - (-\lap)^s_p u(x) = q \xi_\lambda^{q-1}(x) w_\lambda (x) ,$$
where $\xi_\lambda (x)$ is valued between $u(x)$ and $u_\lambda(x)$. Hence at points $x \in \Omega_\lambda$ where
$w_\lambda(x) \leq 0$, we have
\begin{equation}
(-\lap)^s_p u_{\lambda}(x) - (-\lap)^s_p u(x) \geq q u^{q-1}(x) w_\lambda (x).
\label{4ineq}
\end{equation}

As usual, we will carry on the proof in two steps. Chose any ray from the origin as the positive $x_1$ direction. First we show that for $\lambda>-1$ but sufficiently close to $-1$, we have
\begin{equation}
w_{\lambda}(x) \geq 0, \;\; \forall \, x \in \Omega_\lambda.
\label{4w}
\end{equation}
This provides the starting point to move the plane. Then in step 2, we move the plane $T_{\lambda}$ toward the right as long as inequality (\ref{4w}) holds to its limiting position. Define
\begin{equation}
\lambda_o = \sup \{\lambda \leq 0 \mid w_{\mu}(x) \geq 0, \, x \in \Omega_{\mu}, \mu \leq \lambda \}.
\label{4.1}
\end{equation}
We will show that $\lambda_o =0$. Since $x_1$ direction is arbitrary, we deduce that the solution $u$ must be radially symmetric and monotone decreasing about the origin. We now show the details in the two steps.

{\em Step 1.} We show that for $\lambda> -1$ and sufficiently closed to $-1$, (\ref{4w}) holds.
Suppose otherwise, then there exists a point $x^o \in \Omega_{\lambda}$, such that
$$ w_{\lambda}(x^o) = \min_{\Omega_\lambda} w_\lambda = \min_{\Sigma_\lambda} w_\lambda < 0.$$
Then
\begin{eqnarray}
& & (-\lap)^s_p u_{\lambda}(x^o) - (-\lap)^s_p u(x^o) \nonumber \\
&=& C_{n, sp} PV \int_{\Sigma_\lambda} \left[ \frac{1}{|x^o-y|^{n+sp}} -  \frac{1}{|x^o-y^\lambda|^{n+sp}}\right] \left[G(u_{\lambda}(x^o) -u_{\lambda}(y)) - G(u(x^o)-u(y))\right]  d y\nonumber \\
&+& C_{n, sp}  \int_{\Sigma_\lambda} \frac{[G(u_{\lambda}(x^o) -u_{\lambda}(y))- G(u(x^o)-u_{\lambda}(y))] + [G(u_{\lambda}(x^o) -u(y)) - G(u(x^o)-u(y))]}{|x^o-y^{\lambda}|^{n+sp}} d y \nonumber\\
&=& C_{n, sp} PV \int_{\Sigma\lambda} \left[ \frac{1}{|x^o-y|^{n+sp}} -  \frac{1}{|x^o-y^\lambda|^{n+sp}}\right] \left[G(u_{\lambda}(x^o) -u_{\lambda}(y)) - G(u(x^o)-u(y))\right]  d y \nonumber\\
&+& C_{n, sp} w_\lambda(x^o)  \int_{\Sigma_\lambda} \frac{G'(\xi(y)) + G'(\eta(y))}{|x^o-y^{\lambda}|^{n+sp}} d y \nonumber \\
& \leq & C_{n, sp} w_\lambda(x^o)  \int_{\Sigma_\lambda} \frac{G'(\xi(y)) + G'(\eta(y))}{|x^o-y^{\lambda}|^{n+sp}} d y = C_{n, sp} w_\lambda(x^o) \, I.
\label{4.10}
\end{eqnarray}
Here
$$t_1(y) < \xi(y) < t_2(y)\;\; \mbox{ and } \;\; t_3(y) < \eta(y) < t_4(y),$$
with
$$ t_1(y) = u_{\lambda}(x^o) -u_{\lambda}(y), \; t_2(y)=u(x^o)-u_{\lambda}(y), \;
t_3(y) = u_{\lambda}(x^o) -u(y), \;  t_4(y)=u(x^o)-u(y).$$

We now estimate $I$ to derive a contradiction with (\ref{4ineq}) when $\lambda$ is sufficiently close to $-1$.

In \cite{CLLg}, under the assumption that
\begin{equation}
G'(t) \geq c_o >0,
\label{4.11}
\end{equation}
one can derive that there exists $c_1>0$, such that
$$ I \geq \frac{c_1}{\delta^{sp}},$$
where $\delta =|\lambda+1|$ is the width of $\Omega_\lambda$ in $x_1$ direction. Hence a contradiction with (\ref{4ineq}) can be deduced for sufficiently small $\delta$. This is known as the {\em Narrow Region Principle}. However, in our case, $G(t)=|t|^{p-2} t$ no longer satisfies (\ref{4.11}), hence we need to conduct a much finer analysis. To this end, we need the following analysis lemma.

\begin{lem}
For $G(t)=|t|^{p-2} t$, it is well-known that by the {\em mean value theorem}, we have
$$ G(t_2) - G(t_1) = G'(\xi) (t_2 -t_1).$$
Then there exists a constant $c_o > 0$, such that
\begin{equation}
|\xi| \geq c_o \max \{ |t_1|, |t_2| \}.
\label{4.12}
\end{equation}
\label{lem4.1}
\end{lem}

The proof of this lemma is quite elementary, while for readers' convenience, we include it in the Appendix.

Now we employ this lemma to estimate $I$ in (\ref{4.10}). Let $D = \Sigma_\lambda \setminus \Omega_\lambda$. Noticing that $u(y)=0$ in $D$, we have, for some positive constant $c_1$ and $c_2$,
\begin{equation}
I \geq \int_{\Sigma_\lambda} \frac{c_1 |t_4(y)|^{p-2}}{|x^o -y^\lambda|^{n+sp}} dy \geq
c_1 \int_D \frac{u^{p-2}(x^o)}{|x^o -y^\lambda|^{n+sp}} dy \geq \frac{c_2 u^{p-2}(x^o)}{\delta^{sp}} .
\label{4.13}
\end{equation}
Here $\delta = (\lambda +1)$ is the width of the region $\Omega_\lambda$ in the $x_1$-direction.

It follows from (\ref{4.10}) and (\ref{4.13}) that
\begin{eqnarray*}
& & (-\lap)^s_p (u_{\lambda}(x^o)) - (-\lap)^s_p (u(x^o)) - q u^{q-1}(x^o) w_\lambda (x^o)\\
&\leq& C_{n,sp} w_\lambda (x^o) \left[ \frac{c_2 u^{p-2}(x^o)}{\delta^{sp}} - q u^{q-1}(x^o) \right] \\
&<& 0,
\end{eqnarray*}
when $\delta$ is sufficiently small and when $q \geq p-1$.  This contradicts (\ref{4ineq}). Therefore (\ref{4w}) must be true for $\lambda$ sufficiently close to $-1$.
\medskip

{\em Step 2.} Step 1 provides a starting point to move the plane $T_\lambda$. Now we move the plane to the right as long as (\ref{4w}) holds to its limiting position. More precisely, define
$$ \lambda_o = \sup \{\lambda \leq 0 \mid w_\mu (x) \geq 0, x \in \Omega_\mu, \mu \leq \lambda \},$$
and we will show that
\begin{equation}
\lambda_o = 0 .
\label{4.14}
\end{equation}

Suppose in the contrary, $\lambda_o < 0$, then obviously $w_{\lambda_o}$ is not identically zero in
$\Omega_{\lambda_o}$, and by the {\em strong maximum principle} (second part of Theorem \ref{thm2.2}), we have
\begin{equation}
w_{\lambda_o}(x) > 0 , \;\; \forall \, x \in \Omega_{\lambda_o} .
\label{4.15}
\end{equation}
To see this, we argue by contradiction. Suppose there is a point $x^o \in \Omega_{\lambda_o}$, such that $w_{\lambda_o} =0$, then from the equation, we derive
\be
(-\lap)^s_p u_{\lambda}(x^o) - (-\lap)^s_p u(x^o) = 0.
\label{4.16}
\ee
On the other hand, from the proof of the second part of Theorem \ref{thm2.2}, we deduce
$$(-\lap)^s_p u_{\lambda}(x^o) - (-\lap)^s_p u(x^o) < 0.$$
A contradiction with (\ref{4.16}). Hence (\ref{4.15}) must be valid.

By the definition of $\lambda_o$, there exists a sequence $0 \geq \lambda_k \searrow \lambda_o$, and $x^k \in \Omega_{\lambda_k}$, such that
\begin{equation}
w_{\lambda_k}(x^k) = \min_{\Sigma_{\lambda_k}} w_{\lambda_k}  < 0, \;\; \mbox{ and } \; \grad w_{\lambda_k}(x^k) = 0 .
\label{4.17}
\end{equation}

There is a subsequence of $\{x^k\}$ that converges to some point $x^o$, and from (\ref{4.17}) and the continuity of $w_\lambda(x)$ and its derivative with respect to both $x$ and $\lambda$, we have
\be w_{\lambda_o}(x^o) \leq 0, \mbox{ hence } x^o \in \partial \Sigma_{\lambda_o}; \;\; \mbox{ and } \; \grad w_{\lambda_o}(x^o) = 0 .
\label{4.18}
\ee

It follows from equation (\ref{mfu1a}),
\be
\frac{1}{\delta_k} \left\{ (-\lap)^s_p u_{\lambda_k} (x^k) - (-\lap)^s_p u(x^k) \right\} = \frac{w_{\lambda_k}(x^k)}{\delta_k}.
\label{4.19}
\ee
(\ref{4.18}) implies that the right hand side of (\ref{4.19}) must go to zero, and
his will contradicts Theorem \ref{thmHopf}. Therefore, we must have
$$ \lambda_o = 0, \;\; \mbox{ and } \;\; w_0(x) \geq 0 \; \forall \, x \in \Omega_0.$$
Since $x_1$ direction can be chosen arbitrarily, we conclude that $u$ is radially symmetric about the origin.

This completes the proof of the theorem.

\section{Symmetry and Monotonicity in $\mathbb{R}^n$.}

In this section, we study
\be
(-\lap)^s_p u(x) = g(u(x)), \;\; x \in \mathbb{R}^n .
\label{meqwsa}
\ee

\begin{thm}
Assume that $u \in C^{1,1}_{loc} \cap L_{sp}$ is a positive solution of (\ref{meqwsa}) with
$\lim_{|x| \ra \infty} u(x) = 0.$
Suppose $g'(s) \leq 0$ for $s>0$ sufficiently small.

Then $u$ must be radially symmetric and monotone decreasing about some point in $\mathbb{R}^n$.
\label{thm5.1}
\end{thm}

{\bf Proof}.

Let $\Sigma_{\lambda}$, $u_{\lambda}(x)$, and $w_{\lambda}(x)$ be defined as in Section 2.

{\em Step 1}. In this step, we show that for $\lambda$ sufficiently negative, it holds
\begin{equation}
w_{\lambda}(x) \geq 0, \;\; \forall \, x \in \Sigma_{\lambda} .
\label{5.1}
\end{equation}

From equation (\ref{meqwsa}), we deduce
\begin{equation}
(-\lap)^s_p u_\lambda (x) - (-\lap)^s_p u(x) = g'(\xi_\lambda (x)) w_\lambda (x) ,
\label{5.2}
\end{equation}
where $\xi_\lambda (x)$ values between $u_\lambda (x)$ and $u(x)$.

Suppose that (\ref{5.1}) is violated, then by the condition $\lim_{|x| \ra \infty} u(x) = 0$, there exists a point $x^o \in \Sigma_\lambda$, such that
$$ w_{\lambda}(x^o) = \min_{\Sigma_\lambda} w_\lambda < 0 .$$
Consequently,
$$ u_\lambda (x^o) \leq \xi_\lambda (x^o) \leq u(x^o) .$$

For sufficiently negative $\lambda$, $u(x^o)$ is small, hence is $\xi_\lambda (x^o)$, and therefore
$g'(\xi_\lambda (x^o)) \leq 0$ due to the condition on $g'(\cdot)$. It follows that
\begin{equation}
(-\lap)^s_p u_\lambda (x^o) - (-\lap)^s_p u(x^o) \geq 0.
\label{5.3}
\end{equation}

While on the other hand, from the proof of Theorem \ref{thm2.2}, we derive
\be
(-\lap)^s_p u_\lambda (x^o) - (-\lap)^s_p u(x^o) < 0.
\label{5.4}
\ee
This contradicts (\ref{5.3}). Therefore (\ref{5.1}) must be true.
\medskip

{\em Step 2.} Inequality (\ref{5.1}) provides a starting point, from which we move the plane $T_\lambda$ toward the right as long as (\ref{5.1}) holds to its limiting position to show that $u$ is symmetric about the limiting plane. More precisely, let
$$ \lambda_o = \sup \{ \lambda \mid w_\mu (x) \geq 0, \; x \in \Sigma_\mu, \; \mu \leq \lambda \}, $$
we show that $u$ is symmetric about the limiting plane $T_{\lambda_o}$, or
\begin{equation}
 w_{\lambda_o}(x) \equiv 0 , \;\; x \in \Sigma_{\lambda_o}.
 \label{5.5}
 \end{equation}
Suppose (\ref{5.5}) is false, then by the {\em strong maximum principle} (see the argument in Step 2 in the proof of Theorem \ref{thm4.1}), we have
$$w_{\lambda_o} (x) > 0 , \;\; \forall \, x \in \Sigma_{\lambda_o} .$$

On the other hand,  by the definition of $\lambda_o$, there exists a sequence $\lambda_k \searrow \lambda_o$, and $x^k \in \Sigma_{\lambda_k}$, such that
\begin{equation}
w_{\lambda_k}(x^k) = \min_{\Sigma_{\lambda_k}} w_{\lambda_k}  < 0, \;\; \mbox{ and } \; \grad w_{\lambda_k}(x^k) = 0 .
\label{5.6}
\end{equation}

Then the assumption
$$ g'(s) \leq 0, \; \mbox{ for } s \mbox{ small }$$
guarantees that there is a subsequence of $\{x^k\}$ that converges to some point $x^o$.
In fact, from
\begin{equation}
(-\lap)^s_p u_{\lambda_k} (x^k) - (-\lap)^s_p u(x^k) = g'(\xi_{\lambda_k} (x^k)) w_{\lambda_k} (x^k) ,
\label{5.7}
\end{equation}
one can see that if $|x^k|$ is sufficiently large, $u(x^k)$ is small, hence $\xi_{\lambda_k} (x^k)$ is small, this implies that $g'(\xi_{\lambda_k} (x^k)) \leq 0$, and consequent by (\ref{5.7}),
$$(-\lap)^s_p u_{\lambda_k} (x^k) - (-\lap)^s_p u(x^k) \geq 0. $$
This contradicts the fact that $x^k$ is a negative minimum of $w_{\lambda_k}$ (see the proof of Theorem \ref{thm2.2}). Therefore $\{x^k\}$ must be bounded.

Now from (\ref{5.6}), we have
$$w_{\lambda_o}(x^o) \leq 0, \mbox{ hence } x^o \in \partial \Sigma_{\lambda_o}; \;\; \mbox{ and } \; \grad w_{\lambda_o}(x^o) = 0 .$$
It follows that
$$ \frac{w_{\lambda_k}(x^k)}{\delta_k} \ra 0, \;\; \mbox{ as } k \ra \infty.$$
This contradicts Theorem \ref{thmHopf} via (\ref{5.7}) and hence completes the proof of Theorem \ref{thm5.1}.

\section{Appendix}

In this section, we prove the analysis lemma.

\begin{lem}
For $G(t)=|t|^{p-2} t$, it is well-known that by the {\em mean value theorem}, we have
$$ G(t_2) - G(t_1) = G'(\xi) (t_2 -t_1).$$
Then there exists a constant $c_o > 0$, such that
\begin{equation}
|\xi| \geq c_o \max \{ |t_1|, |t_2| \}.
\label{6.12}
\end{equation}
\label{lem6.1}
\end{lem}

{\bf Proof.}

Without loss of generality, we may assume that $|t_2| \geq |t_1|$.

{\em Case i)} $|t_1| \geq \frac{|t_2|}{2}.$

If $t_1$ and $t_2$ are of the same sign, then
$\xi$ is between $t_1$ and $t_2$, we are done. We now assume that $t_1$ and $t_2$ are of opposite signs, hence $G(t_1)$ and $G(t_2)$ are of opposite signs. It follows that
\begin{eqnarray*}
2(p-1)|\xi|^{p-2}|t_2| \geq G'(\xi)|t_2-t_1| =|G(t_2)-G(t_1)| \geq |G(t_2)| =|t_2|^{p-1}.
\end{eqnarray*}
This implies (\ref{6.12}).

{\em Case ii)} $|t_1| \leq \frac{|t_2|}{2}.$

Then
$$2(p-1)|\xi|^{p-2}|t_2| \geq |G(t_2)-G(t_1)| \geq |G(t_2)|-|G(t_1)| \geq c_1 |G(t_2)|.$$
This also implies (\ref{6.12}), and hence completes the proof of the lemma.

\bigskip

{\em Authors' Addresses and E-mails:}
\medskip

Wenxiong Chen

%Department of Applied Mathematics

%Northwestern Polytechnical University

%Xi'an,   China  and

Department of Mathematical Sciences

Yeshiva University

New York, NY, 10033 USA

wchen@yu.edu
\medskip

Congming Li

School of Mathematical Sciences

Institute of Natural Sciences, and MOE-LSC

Shanghai Jiao Tong University

Shanghai, China, and

Department of Applied Mathematics

University of Colorado,

Boulder CO USA

congmingli@gmail.com

\end{document}